\documentclass[11pt]{article}
\usepackage{amssymb}
\usepackage{epsfig}

\def\be{\begin{equation}}
\def\ee{\end{equation}}
\def\bea{\begin{eqnarray}}
\def\eea{\end{eqnarray}}
\def\bes{\begin{eqnarray*}}
\def\ees{\end{eqnarray*}}

\def\nn{\nonumber}
\def\<{\langle}
\def\>{\rangle}
\def\lb{\label}
\def\bs{\setminus}

\def\R{{\bf R}}

\def\Z{{\bf Z}}
\def\N{{\bf N}}

\def\Q{{\bf Q}}
\def\T{{\bf T}}

\def\aa{{\alpha}}

\def\ga{{\gamma}}

\def\ep{{\epsilon}}
\def\lm{{\lambda}}
\def\Lm{{\Lambda}}

\def\Sg{{\Sigma}}
\def\vf{{\varphi}}

\def\H{{\cal H}}
\def\T{{\cal T}}

\def\Nn{{\cal N}}

\def\mul{{\rm mul}}

\def\crit{{\rm crit}}

\def\Sp{{\rm Sp}}

\def\wtd#1{\widetilde{#1}}

\title{The existence of two closed characteristics on every compact star-shaped hypersurface in ${\bf R}^4$}
\author{Hui Liu$^{1}$\thanks{Partially supported by China Postdoctoral Science Foundation No.2013M540512..
E-mail:huiliu@ustc.edu.cn } \qquad  and  \qquad
Yiming Long$^{2}$\thanks{Partially supported by NSFC (No. 11131004), MCME, LPMC of MOE
of China, Nankai University and BCMIIS of Capital Normal University. E-mail: longym@nankai.edu.cn }  \\ \\
$^{1}$ School of Mathematical Sciences, University of Science and Technology of China, \\
Hefei, Anhui 230026, People's Republic of China\\
$^{2}$ Chern Institute of Mathematics and LPMC, Nankai University, \\Tianjin 300071, People's Republic of China\\}
\date{}

\begin{document}

\maketitle

\begin{abstract}
{\it  Recently, Cristofaro-Gardiner and Hutchings proved that there exist at least two closed
characteristics on every compact star-shaped hypersuface in ${\bf R}^4$. Then Ginzburg, Hein,
Hryniewicz, and Macarini gave this result a second proof. In this paper, we give it a third proof
by using index iteration theory, resonance identities of closed characteristics and a remarkable
theorem of Ginzburg et al.}
\end{abstract}



\renewcommand{\theequation}{\thesection.\arabic{equation}}
\renewcommand{\thefigure}{\thesection.\arabic{figure}}

\setcounter{equation}{0}
\section{Introduction}

Let $\Sigma$ be a $C^3$ compact hypersurface in $\R^{2n}$ strictly star-shaped with respect to the origin, i.e.,
the tangent hyperplane at any $x\in\Sigma$ does not intersect the origin. We denote the set of all such
hypersurfaces by $\H_{st}(2n)$, and denote by $\H_{con}(2n)$ the subset of $\H_{st}(2n)$ which consists of all
strictly convex hypersurfaces. We consider closed characteristics $(\tau, y)$ on $\Sigma$, which are solutions
of the following problem
\be
\left\{\matrix{\dot{y}=JN_\Sigma(y), \cr
               y(\tau)=y(0), \cr }\right. \lb{1.1}\ee
where $J=\left(\begin{array}{cc}
                0 & -I_{n} \\
                I_{n} & 0 \\ \end{array}\right)$, $I_n$ is the identity matrix in $\R^n$, $\tau>0$, and
$N_\Sigma(y)$ is the outward normal vector of $\Sigma$ at $y$ normalized by the condition $N_\Sigma(y)\cdot y=1$.
Here $a\cdot b$ denotes the standard inner product of $a, b\in\R^{2n}$. A closed characteristic $(\tau, y)$ is
{\it prime}, if $\tau$ is the minimal period of $y$. Two closed characteristics $(\tau, y)$ and $(\sigma, z)$
are {\it geometrically distinct}, if $y(\R)\not= z(\R)$. We denote by $\T(\Sg)$ the set of all geometrically
distinct closed characteristics on $\Sg$.

Fix a constant $\alpha$ satisfying $1<\alpha<2$ and define the Hamiltonian function $H_\alpha:\R^{2n}\to [0,+\infty)$ by
\be  H_\alpha(x) = j(x)^{\alpha}, \qquad \forall x\in\R^{2n},  \lb{1.2}\ee
where $j$ is the gauge function of $\Sigma$, i.e., $j(x)=\lm$ if $x=\lm y$ for some $\lm>0$ and $y\in\Sigma$ when
$x\in\R^{2n}\bs\{0\}$, and $j(0)=0$. Then $H_{\alpha}\in C^1(\R^{2n},\R)\cap C^3(\R^{2n}\bs\{0\},\R)$ and
$\Sigma=H_{\alpha}^{-1}(1)$. It is well-known that the problem (\ref{1.1}) is equivalent to the following given energy
problem of the Hamiltonian system
\bea \cases{\dot{y}(t) &$=JH_{\alpha}^{\prime}(y(t)), \quad H_{\alpha}(y(t)) =1, \qquad\forall t\in\R,$ \cr
             y(\tau) &$= y(0).$ \cr} \lb{1.3}\eea
Denote by $\T(\Sigma,\alpha)$ the set of all geometrically distinct solutions $(\tau,y)$ of the problem (\ref{1.3}).
Note that elements in $\T(\Sigma)$ and $\T(\Sigma,\alpha)$ are one to one correspondent to each other.

Let $(\tau,y)\in \T(\Sigma,\alpha)$. The fundamental solution $\gamma_y:[0,\tau]\to\Sp(2n)$ with $\gamma_y(0)=I_{2n}$
of the linearized Hamiltonian system
\bea \dot{\xi}(t) = JH_{\alpha}''(y(t))\xi(t), \qquad \forall t\in\R,  \nn\eea
is called the {\it associated symplectic path} of $(\tau,y)$. The eigenvalues of $\gamma_y(\tau)$ are called
{\it Floquet multipliers} of $(\tau,y)$. By Proposition 1.6.13 of \cite{Eke1}, the Floquet multipliers with their
multiplicities of $(\tau,y)\in\mathcal{T}(\Sigma,\alpha)$ do not depend on the
particular choice of the Hamiltonian function in (\ref{1.2}). A closed characteristic $(\tau,y)$ is
{\it non-degenerate}, if $1$ is a Floquet multiplier of $y$ of precisely algebraic multiplicity $2$, and
is {\it hyperbolic}, if $1$ is a double Floquet multiplier of it and all the other Floquet
multipliers are not on ${\bf U}=\{z\in {\bf C}\mid |z|=1\}$, i.e., the unit circle in the complex plane,
and is {\it elliptic}, if all the Floquet multipliers of $y$ are on ${\bf U}$.

The study on closed characteristics in the global sense started in 1978, when the existence of at least one
closed characteristic was first established on any $\Sg\in\H_{st}(2n)$ by P. Rabinowitz in \cite{Rab1}
and on any $\Sg\in\H_{con}(2n)$ by A. Weinstein in \cite{Wei1} independently. Since then the existence of
multiple closed characteristics on $\Sg\in\H_{con}(2n)$ has been deeply studied by many mathematicians, for
example, studies in \cite{EkL1}, \cite{EkH1}, \cite{Szu1}, \cite{HWZ1}, \cite{LoZ1}, \cite{WHL1}, \cite{Wan1} for
convex hypersurfaces. In \cite{HuL1} of 2002, X. Hu and Y. Long proved $\;^{\#}\T(\Sg)\ge 2$ for every $\Sg\in\H_{st}(2n)$
with $n\ge 2$ provided all the closed characteristics on $\Sg$ and their iterates are non-degenerate.
Note that in October of 2011, the authors and W. Wang established the resonance identities for closed characteristics
on compact star-shaped hypersurface in $\R^{2n}$ (cf. \cite{LLW1}), and the authors proved that any symmetric compact
star-shaped hypersurface in $\R^4$ has at least two closed characteristics regardless of whether the resulting contact
form is non-degenerate or not, which was announced in October 2011 at a conference in Princeton (cf. \cite{Lon3}).
Shortly after that in February of 2012, $\;^{\#}\T(\Sg)\ge 2$ was proved for every $\Sg\in \H_{st}(4)$ by
D. Cristofaro-Gardiner and M. Hutchings in \cite{CGH1}. More recently, V. Ginzburg, D. Hein, U. Hryniewicz, and L.
Macarini proved a remarkable theorem in \cite{GHHM} which asserts that the existence of one simple closed Reeb orbit of
a particular type (a symplectically degenerate maximum) forces the Reeb flow to have infinitely many periodic orbits.
Then by using this theorem they also proved $\;^{\#}\T(\Sg)\ge 2$ for every $\Sg\in \H_{st}(4)$.

In this paper, we give a third proof for this result.

{\bf Theorem 1.1.} {\it $\,^{\#}\T(\Sg)\geq 2$ holds for every $\Sg\in \H_{st}(4)$. }

{\bf Remark 1.2.} The resonance identities for closed characteristics on compact star-shaped hypersurface
in $\R^{2n}$ established in \cite{LLW1} serve as a useful tool to study the multiplicity
and stability of closed characteristics on star-shaped hypersurfaces. Besides these mean index identities, the other
main ingredients in the proof of Theorems 1.1 are: Morse inequality, the index iteration theory developed by Y. Long
and his coworkers (cf. \cite{Lon2}), specially the precise iteration formulae of the Maslov-type index theory for
any symplectic path established by Y. Long in \cite{Lon1}, and the remarkable Theorem 1.2 of \cite{GHHM} about the
symplectically degenerate maximum. In Section 2, we briefly review the equivariant Morse theory and the resonance
identities for closed characteristics on compact star-shaped hypersurfaces in $\R^{2n}$ developed in
\cite{LLW1}. In Section 3, we give the proof of Theorem 1.1. 

In this paper, let $\N$, $\N_0$, $\Z$, $\Q$, $\R$, and $\R^+$ denote the sets of natural integers,
non-negative integers, integers, rational numbers, real numbers, and positive real numbers respectively.
We define the functions $[a]=\max{\{k\in {\bf Z}\mid k\leq a\}}$, $\{a\}=a-[a]$ , and
$E(a)=\min{\{k\in{\bf Z}\mid k\geq a\}}$. Denote by $a\cdot b$ and $|a|$ the standard
inner product and norm in $\R^{2n}$. Denote by $\langle\cdot,\cdot\rangle$ and $\|\cdot\|$
the standard $L^2$ inner product and $L^2$ norm. For an $S^1$-space $X$, we denote by
$X_{S^1}$ the homotopy quotient of $X$ by $S^1$, i.e., $X_{S^1}=S^\infty\times_{S^1}X$,
where $S^\infty$ is the unit sphere in an infinite dimensional {\it complex} Hilbert space.
In this paper we use only $\Q$ coefficients for all homological modules. By $t\to a^+$, we
mean $t>a$ and $t\to a$. For a $\Z_m$-space pair $(A, B)$, let
$H_{\ast}(A, B)^{\pm\Z_m} = \{\sigma\in H_{\ast}(A, B)\,|\,L_{\ast}\sigma=\pm \sigma\}$,
where $L$ is a generator of the $\Z_m$-action.

\setcounter{equation}{0}
\section{Equivariant Morse theory for closed characteristics on compact star-shaped hypersurfaces}
In the rest of this paper, we fix a $\Sg\in\H_{st}(2n)$ and assume the following condition on $\T(\Sg)$:

\noindent (F) {\bf There exist only finitely many geometrically distinct prime closed characteristics\\
$\qquad\qquad \{(\tau_j, y_j)\}_{1\le j\le k}$ on $\Sigma$. }

In this section, we briefly review the equivariant Morse theory, especially the resonance identities
for closed characteristics on $\Sg\in\H_{st}(2n)$ established in \cite{LLW1} which will be needed in the
proof in Section 3. All the details of proofs in this section can be found in \cite{LLW1}.

Let $\hat{\tau}=\inf_{1\leq j\leq k}{\tau_j}$ and $T$ be a fixed positive constant. Then by Section 2 of
\cite{LLW1}, for any $a>\frac{\hat{\tau}}{T}$, we can construct a function $\varphi_a\in C^{\infty}({\bf R}, {\bf R}^+)$
which has $0$ as its unique critical point in $[0, +\infty)$. Moreover, $\frac{\varphi^{\prime}(t)}{t}$ is strictly
decreasing for $t>0$ together with $\varphi(0)=0=\varphi^{\prime}(0)$ and
$\varphi^{\prime\prime}(0)=1=\lim_{t\rightarrow 0^+}\frac{\varphi^{\prime}(t)}{t}$. More precisely, we
define $\varphi_a$ and the Hamiltonian function $\wtd{H}_a(x)=a\vf_a(j(x))$ via Lemma 2.2 and Lemma 2.4
in \cite{LLW1}. The precise dependence of $\varphi_a$ on $a$ is explained in Remark 2.3 of \cite{LLW1}.

For technical reasons we want to further modify the Hamiltonian, we define the new Hamiltonian
function $H_a$ via Proposition 2.5 of \cite{LLW1} and consider the fixed period problem
\be  \dot{x}(t)=JH_a^\prime(x(t)),\quad x(0)=x(T).  \lb{2.1}\ee
Then $H_a\in C^{3}({\bf R}^{2n} \setminus\{0\},{\bf R})\cap C^{1}({\bf R}^{2n},{\bf R})$.
Solutions of (\ref{2.1}) are $x\equiv 0$ and $x=\rho y(\tau t/T)$ with
$\frac{\vf_a^\prime(\rho)}{\rho}=\frac{\tau}{aT}$, where $(\tau, y)$ is a solution of (\ref{1.1}). In particular,
non-zero solutions of (\ref{2.1}) are in one to one correspondence with solutions of (\ref{1.1}) with period
$\tau<aT$.

For any $a>\frac{\hat{\tau}}{T}$, we can choose some large constant $K=K(a)$ such that
\be H_{a,K}(x) = H_a(x)+\frac{1}{2}K|x|^2   \lb{2.2}\ee
is a strictly convex function, that is,
\be (\nabla H_{a, K}(x)-\nabla H_{a, K}(y), x-y) \geq \frac{\ep}{2}|x-y|^2,  \lb{2.3}\ee
for all $x, y\in {\bf R}^{2n}$, and some positive $\ep$. Let $H_{a,K}^*$ be the Fenchel dual of $H_{a,K}$
defined by
$$  H_{a,K}^\ast (y) = \sup\{x\cdot y-H_{a,K}(x)\;|\; x\in \R^{2n}\}.   $$
The dual action functional on $X=W^{1, 2}({\bf R}/{T {\bf Z}}, {\bf R}^{2n})$ is defined by
\be F_{a,K}(x) = \int_0^T{\left[\frac{1}{2}(J\dot{x}-K x,x)+H_{a,K}^*(-J\dot{x}+K x)\right]dt}. \lb{2.4}\ee
Then $F_{a,K}\in C^{1,1}(X, \R)$ and for $KT\not\in 2\pi{\bf Z}$, $F_{a,K}$ satisfies the
Palais-Smale condition and $x$ is a critical point of $F_{a, K}$ if and only if it is a solution of (\ref{2.1}). Moreover,
$F_{a, K}(x_a)<0$ and it is independent of $K$ for every critical point $x_a\neq 0$ of $F_{a, K}$.

When $KT\notin 2\pi{\bf Z}$, the map $x\mapsto -J\dot{x}+Kx$ is a Hilbert space isomorphism between
$X=W^{1, 2}({\bf R}/({T {\bf Z}}); {\bf R}^{2n})$ and $E=L^{2}({\bf R}/(T {\bf Z}),{\bf R}^{2n})$. We denote its inverse
by $M_K$ and the functional
\be \Psi_{a,K}(u)=\int_0^T{\left[-\frac{1}{2}(M_{K}u, u)+H_{a,K}^*(u)\right]dt}, \qquad \forall\,u\in E. \lb{2.5}\ee
Then $x\in X$ is a critical point of $F_{a,K}$ if and only if $u=-J\dot{x}+Kx$ is a critical point of $\Psi_{a, K}$.

Suppose $u$ is a nonzero critical point of $\Psi_{a, K}$.
Then the formal Hessian of $\Psi_{a, K}$ at $u$ is defined by
\be Q_{a,K}(v)=\int_0^T(-M_K v\cdot v+H_{a,K}^{*\prime\prime}(u)v\cdot v)dt,  \lb{2.6}\ee
which defines an orthogonal splitting $E=E_-\oplus E_0\oplus E_+$ of $E$ into negative, zero and positive subspaces.
The index and nullity of $u$ are defined by $i_K(u)=\dim E_-$ and $\nu_K(u)=\dim E_0$ respectively.
Similarly, we define the index and nullity of $x=M_Ku$ for $F_{a, K}$, we denote them by $i_K(x)$ and
$\nu_K(x)$. Then we have
\be  i_K(u)=i_K(x),\quad \nu_K(u)=\nu_K(x),  \lb{2.7}\ee
which follow from the definitions (\ref{2.4}) and (\ref{2.5}). The following important formula was proved in
Lemma 6.4 of \cite{Vit2}:
\be  i_K(x) = 2n([KT/{2\pi}]+1)+i^v(x) \equiv d(K)+i^v(x),   \lb{2.8}\ee
where the index $i^v(x)$ does not depend on K, but only on $H_a$.

By the proof of Proposition 2 of \cite{Vit1}, we have that $v\in E$ belongs to the null space of $Q_{a, K}$
if and only if $z=M_K v$ is a solution of the linearized system
\be  \dot{z}(t) = JH_a''(x(t))z(t).  \lb{2.9}\ee
Thus the nullity in (\ref{2.7}) is independent of $K$, which we denote by $\nu^v(x)\equiv \nu_K(u)= \nu_K(x)$.

By Proposition 2.11 of \cite{LLW1}, the index $i^v(x)$ and nullity $\nu^v(x)$ coincide with those defined for
the Hamiltonian $H(x)=j(x)^\alpha$ for all $x\in\R^{2n}$ and some $\aa\in (1,2)$. Especially
$1\le \nu^v(x_b)\le 2n-1$ always holds.

We have a natural $S^1$-action on $X$ or $E$ defined by
$$  \theta\cdot u(t)=u(\theta+t),\quad\forall\, \theta\in S^1, \, t\in\R.  $$
Clearly both of $F_{a, K}$ and $\Psi_{a, K}$ are $S^1$-invariant. For any $\kappa\in\R$, we denote by
\bea
\Lambda_{a, K}^\kappa &=& \{u\in L^{2}({\bf R}/({T {\bf Z}}); {\bf R}^{2n})\;|\;\Psi_{a,K}(u)\le\kappa\}  \nn\\
X_{a, K}^\kappa &=& \{x\in W^{1, 2}({\bf R}/(T {\bf Z}),{\bf R}^{2n})\;|\;F_{a, K}(x)\le\kappa\}.  \nn\eea
For a critical point $u$ of $\Psi_{a, K}$ and the corresponding $x=M_K u$ of $F_{a, K}$, let
\bea
\Lm_{a,K}(u) &=& \Lm_{a,K}^{\Psi_{a, K}(u)}
   = \{w\in L^{2}(\R/(T\Z), \R^{2n}) \;|\; \Psi_{a, K}(w)\le\Psi_{a,K}(u)\},  \nn\\
X_{a,K}(x) &=& X_{a,K}^{F_{a,K}(x)} = \{y\in W^{1, 2}(\R/(T\Z), \R^{2n}) \;|\; F_{a,K}(y)\le F_{a,K}(x)\}. \nn\eea
Clearly, both sets are $S^1$-invariant. Denote by $\crit(\Psi_{a, K})$ the set of critical points of $\Psi_{a, K}$.
Because $\Psi_{a,K}$ is $S^1$-invariant, $S^1\cdot u$ becomes a critical orbit if $u\in \crit(\Psi_{a, K})$.
Note that by the condition (F), the number of critical orbits of $\Psi_{a, K}$
is finite. Hence as usual we can make the following definition.

{\bf Definition 2.1.} {\it Suppose $u$ is a nonzero critical point of $\Psi_{a, K}$, and $\Nn$ is an $S^1$-invariant
open neighborhood of $S^1\cdot u$ such that $\crit(\Psi_{a,K})\cap (\Lm_{a,K}(u)\cap \Nn) = S^1\cdot u$.
Then the $S^1$-critical module of $S^1\cdot u$ is defined by
$$ C_{S^1,\; q}(\Psi_{a, K}, \;S^1\cdot u)
=H_{q}((\Lambda_{a, K}(u)\cap\Nn)_{S^1},\; ((\Lambda_{a,K}(u)\setminus S^1\cdot u)\cap\Nn)_{S^1}). $$
Similarly, we define the $S^1$-critical module $C_{S^1,\; q}(F_{a, K}, \;S^1\cdot x)$ of $S^1\cdot x$
for $F_{a, K}$.}

We have the following for critical modules.

{\bf Proposition 2.2.} (cf. Proposition 3.2 of \cite{LLW1}) {\it Let $(\tau,y)$ be a closed characteristic on
$\Sigma$. For any $\frac{\tau}{T}<a_1<a_2<+\infty$, let $K$ be a fixed sufficiently large real number so that
(\ref{2.3}) holds for all $a\in [a_1, a_2]$. Then the critical module $C_{S^1,\; q}(F_{a, K}, \;S^1\cdot x)$ is
independent of the choice of $H_a$ in the sense that if $x_i$ is a solution of (\ref{2.1}) with Hamiltonian
function $H_{a_i}(x)$ with $i=1$ and $2$ respectively such that both $x_1$ and $x_2$ correspond to the same
closed characteristic $(\tau,y)$ on $\Sigma$, then we have}
$$ C_{S^1,\;q}(F_{a_1,K},\;S^1\cdot {x}_1) \cong C_{S^1,\;q}(F_{a_2,K},\;S^1\cdot{x}_2),
\qquad \forall\,q\in \Z. $$ 

Now we fix an $a>\frac{\tau}{T}$, and write $F_K$ and $H$ for $F_{a, K}$ and $H_a$
respectively. We suppose also that $K\in {\bf R}$ satisfy (\ref{2.3}).
Then the critical points of $F_{K}$ which are the solutions of (\ref{2.1}) are the same for any
$K$ satisfying that $K\notin \frac{2\pi}{T}{\Z}$.

{\bf Proposition 2.3.} (cf. Theorem 3.3 of \cite{LLW1}) {\it Suppose $\bar{x}$ is a nonzero critical
point of $F_{K}$. Then the $S^1$-critical module $C_{S^1,d(K)+l}(F_{K},S^1\cdot\bar{x})$ is independent of the
choice of $K$ for $KT \notin 2\pi\Z$, i.e.,
$$ C_{S^1,d(K)+l}(F_{K},S^1\cdot\bar{x}) \cong C_{S^1,d(K^\prime)+l}(F_{K^\prime},S^1\cdot\bar{x}), $$
where $KT$, $K^\prime T \notin 2\pi\Z$, $l\in \Z$, and both $K$ and $K^\prime$ satisfy (\ref{2.3}).}

For the fixed $a>\frac{\tau}{T}$ we let $u_K\neq 0$ be a critical point of $\Psi_{a, K}$ with
multiplicity $\mul(u_K)=m$, that is, $u_K$ corresponds to a closed characteristic $(\tau, y)\subset\Sigma$
with $(\tau, y)$ being the $m$-th iterate of some prime closed characteristic. Precisely, we have
$u_K=-J\dot x+Kx$ with $x$ being a solution of (\ref{2.1}) and $x=\rho y(\frac{\tau t}{T})$ with
$\frac{\vf_a^\prime(\rho)}{\rho}=\frac{\tau}{aT}$.
Moreover, $(\tau, y)$ is a closed characteristic on $\Sigma$ with minimal period $\frac{\tau}{m}$. Hence the isotropy
group of $(\tau, y)$ satisfies $\{\theta\in S^1\;|\;\theta\cdot u_K=u_K\}=\Z_m$ and the orbit of $u_K$ satisfies
$S^1\cdot u_K\cong S^1/\Z_m\cong S^1$. Since $\Psi_{a, K}$ is not $C^2$ on $E=L^{2}({\bf R}/(T{\bf Z}),{\bf R}^{2n})$,
we need to use a finite-dimensional approximation introduced by Viterbo in order to apply Morse theory.
More precisely, by Lemma 2.10 of \cite{LLW1}, we construct a finite dimensional $S^1$-invariant subspace $G$ of
$L^{2}({\bf R}/(T{\bf Z}); {\bf R}^{2n})$ and a functional $\psi_{a,K}$ on $G$. Moreover, $\Psi_{a, K}$
and $\psi_{a,K}$ have the same critical points, $\psi_{a,K}$ is $C^{2}$ in a small tubular neighborhood of
the critical orbit $S^1\cdot g_K$, where $g_K$ is the critical point of $\psi_{a,K}$ corresponding to
$u_K$, and the Morse index and nullity of its critical points coincide with those of the corresponding critical points
of $\Psi_{a, K}$, and then the isotropy group satisfies $\{\theta\in S^1\;|\;\theta\cdot g_K=g_K\}=\Z_m$.
Let $p: N(S^1\cdot g_K)\rightarrow S^1\cdot g_K$ be the normal bundle of $S^1\cdot g_K$ in $G$
and let $p^{-1}(\theta\cdot g_K)=N(\theta\cdot g_K)$ be the fibre over $\theta\cdot g_K$, where $\theta\in S^1$. Let
$DN(S^1\cdot g_K)$ be the $\varrho$ disk bundle of $N(S^1\cdot g_K)$ for some $\varrho>0$ sufficiently small, i.e.,
$DN(S^1\cdot g_K)=\{\xi\in N(S^1\cdot g_K)\;| \; \|\xi\|<\varrho\}$ which is identified by the exponential map with a
subset of $G$, and let $DN(\theta\cdot g_K)=p^{-1}(\theta\cdot g_K)\cap DN(S^1\cdot g_K)$ be the disk over
$\theta\cdot g_K$. Clearly, $DN(\theta\cdot g_K)$ is $\Z_m$-invariant and we have
$DN(S^1\cdot g_K)=DN(g_K)\times_{\Z_m}S^1$ where the $\Z_m$ action is given by
$$(\theta, v, t)\in \Z_m\times DN(g_K)\times S^1\mapsto (\theta\cdot v, \;\theta^{-1}t)\in DN(g_K)\times S^1. $$
Hence for an $S^1$ invariant subset $\Gamma$ of $DN(S^1\cdot g_K)$, we have
$\Gamma/S^1=(\Gamma_{g_K}\times_{\Z_m}S^1)/S^1=\Gamma_{g_K}/\Z_m$, where $\Gamma_{g_K}=\Gamma\cap DN(g_K)$.
Then we have
$$ C_{S^1,\;\ast}(\Psi_{a,K}, \;S^1\cdot u_K)\cong H_\ast((\wtd{\Lambda}_{a,K}(g_K)\cap DN(g_K)),\;
    ((\wtd{\Lambda}_{a,K}(g_K)\setminus\{g_K\})\cap DN(g_K)))^{\Z_m},   $$
where $\wtd{\Lambda}_{a,K}(g_K)=\{g\in G \;|\; \psi_{a, K}(g)\le\psi_{a, K}(g_K)\}$.

For any $p\in\N$ satisfying $p\tau<aT$, we choose $K$
such that $pK\notin \frac{2\pi}{T}\Z$, then the $p$th iteration $u_{pK}^p$ of $u_K$ is given by $-J\dot x^p+pKx^p$,
where $x^p$ is the unique solution of (\ref{2.1}) corresponding to $(p\tau, y)$ and is a critical point of $F_{a, pK}$,
that is, $u_{pK}^p$ is the critical point of $\Psi_{a, pK}$ corresponding to $x^p$.
By computation, $u_{pK}^p(t)=p^{\frac{\alpha-1}{\alpha-2}}u_K(pt)$.
We define the $p$th iteration $\phi^p$ on $L^{2}(\R/(T\Z); {\bf R}^{2n})$ by
$$ \phi^p: v_K(t)\mapsto v^p_{pK}(t)\equiv p^\frac{\alpha-1}{\alpha-2}v_K(pt), \forall p\in\N. $$
Denote by $g_{pK}^p\equiv \phi^p(g_K)$ the critical point of $\psi_{a,pK}$ corresponding to $u_{pK}^p$.

Now we can apply the results of D. Gromoll and W. Meyer \cite{GrM1} to the manifold $DN(g_{pK}^p)$ with $g_{pK}^p$
as the unique critical point. Then $\mul(g_{pK}^p)=pm$ is the multiplicity of $g_{pK}^p$ and the isotropy group
$\Z_{pm}\subseteq S^1$ of $g_{pK}^p$ acts on $DN(g_{pK}^p)$ by isometries. According to Lemma 1 of \cite{GrM1}, we
have a $\Z_{pm}$-invariant decomposition of $T_{g_{pK}^p}(DN(g_{pK}^p))$:
$$  T_{g_{pK}^p}(DN(g_{pK}^p)) = V^+\oplus V^-\oplus V^0 = \{(x_+, x_-, x_0)\},  $$
with $\dim V^-=i(g_{pK}^p)=i_{pK}(u_{pK}^p)$, $\dim V^0=\nu(g_{pK}^p)-1=\nu_{pK}(u_{pK}^p)-1$,
and a $\Z_{pm}$-invariant neighborhood $B=B_+\times B_-\times B_0$ for $0$ in $T_{g_{pK}^p}(DN(g_{pK}^p))$ together
with two $\Z_{pm}$-invariant diffeomorphisms
$$ \Phi : B=B_+\times B_-\times B_0\to \Phi(B_+\times B_-\times B_0)\subset DN(g_{pK}^p),  $$
and
$$ \eta : B_0\to W(g_{pK}^p)\equiv\eta(B_0)\subset DN(g_{pK}^p), $$
and $\Phi(0)=\eta(0)=g_{pK}^p$, such that
$$ \psi_{a,pK}\circ\Phi(x_+,x_-,x_0)=|x_+|^2 - |x_-|^2 + \psi_{a,pK}\circ\eta(x_0), $$
with $d(\psi_{a, pK}\circ \eta)(0)=d^2(\psi_{a, pK}\circ\eta)(0)=0$. As usual, we call $W(g_{pK}^p)$ a local
characteristic manifold, and $U(g_{pK}^p)=B_-$ a local negative disk at $g_{pK}^p$. By the proof of Lemma 1 of
\cite{GrM1}, $W(g_{pK}^p)$ and $U(g_{pK}^p)$ are $\Z_{pm}$-invariant. Then we have
\bea
&& H_\ast(\wtd{\Lambda}_{a, pK}(g_{pK}^p)\cap DN(g_{pK}^p),\;
     (\wtd{\Lambda}_{a, pK}(g_{pK}^p)\setminus\{g_{pK}^p\})\cap DN(g_{pK}^p)) \nn\\
&&\quad = \bigoplus_{q\in\Z}H_q (U(g_{pK}^p),U(g_{pK}^p)\setminus\{g_{pK}^p\}) \nn\\
&&\quad\qquad \otimes H_{\ast-q}(W(g_{pK}^p)\cap \wtd{\Lambda}_{a,pK}(g_{pK}^p),
       (W(g_{pK}^p)\setminus\{g_{pK}^p\})\cap \wtd{\Lambda}_{a,pK}(g_{pK}^p)), \nn\eea
where
$$  H_q(U(g_{pK}^p),U(g_{pK}^p)\setminus\{g_{pK}^p\} )
    =\left\{ \begin{array}{ll}\Q, & {\rm if\;}q=i_{pK}(u_{pK}^p),\\
0, & {\rm otherwise}.\end{array}\right.  $$
Now we have the following proposition.

{\bf Proposition 2.4.}(cf. Proposition 4.2 of \cite{LLW1})
{\it For any $p\in\N$, we choose $K>0$ such that $pK\notin \frac{2\pi}{T}\Z$. Let $u_K\neq 0$
be a critical point of $\Psi_{a, K}$ with $\mul(u_K)=1$, $u_K=-J\dot x+Kx$ with $x$ being a critical point of
$F_{a, K}$. Then for all $q\in\Z$, we have
\bea
&& C_{S^1,\; q}(\Psi_{a,pK},\;S^1\cdot u_{pK}^p) \nn\\
&&\quad\cong \left(\frac{}{}H_{q-i_{pK}(u_{pK}^p)}(W(g_{pK}^p)\cap \wtd{\Lambda}_{a,pK}(g_{pK}^p),(W(g_{pK}^p)
                 \setminus\{g_{pK}^p\})\cap \wtd{\Lambda}_{a,pK}(g_{pK}^p))\right)^{\beta(x^p)\Z_p}, \nn\eea
where $\beta(x^p)=(-1)^{i_{pK}(u_{pK}^p)-i_K(u_K)}=(-1)^{i^v(x^p)-i^v(x)}$. Thus
$$ C_{S^1,\; q}(\Psi_{a,pK},\;S^1\cdot u_{pK}^p)=0, \qquad {\it if}\quad
      q<i_{pK}(u_{pK}^p)\quad {\it or}\quad q>i_{pK}(u_{pK}^p)+\nu_{pK}(u_{pK}^p)-1. $$
In particular, if $u_{pK}^p$ is non-degenerate, i.e., $\nu_{pK}(u_{pK}^p)=1$, then}
\be C_{S^1,\; q}(\Psi_{a,pK},\;S^1\cdot u_{pK}^p)
    = \left\{\begin{array}{ll}\Q, & {\rm if\;}q=i_{pK}(u_{pK}^p)\;{\rm and\;}\beta(x^p)=1, \\
                      0, & {\rm otherwise}.\end{array}\right. \lb{2.10}\ee
We need the following

{\bf Definition 2.5.} {\it For any $p\in\N$, we choose $K$ such that $pK\notin \frac{2\pi}{T}\Z$. Let $u_K\neq 0$ be
a critical point of $\Psi_{a,K}$ with $\mul(u_K)=1$, $u_K=-J\dot x+Kx$ with $x$ being a critical point of $F_{a, K}$.
Then for all $l\in\Z$, let
\bea
k_{l,\pm 1}(u_{pK}^p) &=& \dim\left(\frac{}{}H_{l}(W(g_{pK}^p)\cap
  \wtd{\Lambda}_{a,pK}(g_{pK}^p),(W(g_{pK}^p)\bs\{g_{pK}^p\})\cap \wtd{\Lambda}_{a,pK}(g_{pK}^p))\right)^{\pm\Z_p},
           \quad  \nn\\
k_l(u_{pK}^p) &=& \dim\left(\frac{}{}H_{l}(W(g_{pK}^p)\cap
  \wtd{\Lambda}_{a,pK}(g_{pK}^p),(W(g_{pK}^p)\bs\{g_{pK}^p\})\cap \wtd{\Lambda}_{a,pK}(g_{pK}^p))\right)^{\beta(x^p)\Z_p}.
           \qquad\quad   \nn\eea
Here $k_l(u_{pK}^p)$'s are called critical type numbers of $u_{pK}^p$.  }

By Proposition 2.3, we obtain that $k_l(u_{pK}^p)$ is independent of the choice of $K$ and denote it by $k_l(x^p)$,
here $k_l(x^p)$'s are called critical type numbers of $x^p$.

We have the following properties for critical type numbers:

{\bf Proposition 2.6.}(cf. Proposition 4.6 of \cite{LLW1})
{\it Let $x\neq 0$ be a critical point of $F_{a,K}$ with $\mul(x)=1$ corresponding to a
critical point $u_K$ of $\Psi_{a, K}$. Then there exists a minimal $K(x)\in \N$ such that
\bea
&& \nu^v(x^{p+K(x)})=\nu^v(x^p),\quad i^v(x^{p+K(x)})-i^v(x^p)\in 2\Z,  \qquad\forall p\in \N,  \nn\\
&& k_l(x^{p+K(x)})=k_l(x^p), \qquad\forall p\in \N,\;l\in\Z. \nn\eea
We call $K(x)$ the minimal period of critical modules of iterations of the functional $F_{a, K}$ at $x$. }

For every closed characteristic $(\tau,y)$ on $\Sigma$, let $aT>\tau$ and choose $\vf_a$ as above.
Determine $\rho$ uniquely by $\frac{\vf_a'(\rho)}{\rho}=\frac{\tau}{aT}$. Let $x=\rho y(\frac{\tau t}{T})$.
Then we define the index $i(\tau,y)$ and nullity $\nu(\tau,y)$ of $(\tau,y)$ by
$$ i(\tau,y)=i^v(x), \qquad \nu(\tau,y)=\nu^v(x). $$
Then the mean index of $(\tau, y)$ is defined by
$$   \hat{i}(\tau,y) = \lim_{m\rightarrow\infty}\frac{i(m\tau, y)}{m}.  $$

Note that by Proposition 2.11 of \cite{LLW1}, the index and nullity are well defined and is independent of the
choice of $a$.

For a closed characteristic $(\tau, y)$ on $\Sigma$, we simply denote by $y^m\equiv(m\tau, y)$
the m-th iteration of $y$ for $m\in\N$. By Proposition 2.3, we can define the critical type numbers $k_l(y^m)$
of $y^m$ to be $k_l(x^m)$, where $x^m$ is the critical point of $F_{a, K}$ corresponding to $y^m$. We also
define $K(y)=K(x)$.

{\bf Remark 2.7.} (cf. Remark 4.10 of \cite{LLW1}) {\it The following holds.

(i) Note that $k_l(y^m)=0$ for $l\notin [0, \nu(y^m)-1]$ and it can take only values $0$ or $1$ when $l=0$
or $l=\nu(y^m)-1$.

(ii) $k_0(y^m)=1$ implies $k_l(y^m)=0$ for $1\le l\le \nu(y^m)-1$.

(iii) $k_{\nu(y^m)-1}(y^m)=1$ implies $k_l(y^m)=0$ for $0\le l\le \nu(y^m)-2$.

(iv) $k_l(y^m)\ge 1$ for some $1\le l\le \nu(y^m)-2$ implies $k_0(y^m)=k_{\nu(y^m)-1}(y^m)=0$.

(v) In particular, only one of the $k_l(y^m)$s for $0\le l\le \nu(y^m)-1$ can be non-zero when $\nu(y^m)\le 3$.}

For a closed characteristic $(\tau, y)$ on $\Sigma$, the average Euler characteristic $\hat\chi(y)$ of $y$ is
defined by
$$ \hat\chi(y)=\frac{1}{K(y)}\sum_{1\le m\le K(y)\atop 0\le l\le 2n-2}(-1)^{i(y^{m})+l}k_l(y^{m}).  $$
$\hat\chi(y)$ is a rational number. In particular, if all $y^m$s are non-degenerate, then by Proposition 2.6 we have
$$ \hat\chi(y)
    = \left\{\begin{array}{ll}(-1)^{i(y)}, & {\rm if\;\;} i(y^2)-i(y)\in 2\Z,  \\
           \frac{(-1)^{i(y)}}{2}, & {\rm otherwise}. \end{array}\right.  $$
We have the following resonance identities for closed characteristics.

{\bf Theorem 2.8.} {\it Suppose that $\Sg\in \H_{st}(2n)$ satisfies
$\,^{\#}\T(\Sg)<+\infty$. Denote all the geometrically distinct prime closed characteristics by
$\{(\tau_j,\; y_j)\}_{1\le j\le k}$. Then the following identities hold
\bea
\sum_{1\le j\le k\atop \hat{i}(y_j)>0}\frac{\hat{\chi}(y_j)}{\hat{i}(y_j)} &=& \frac{1}{2},  \nn\\
\sum_{1\le j\le k\atop \hat{i}(y_j)<0}\frac{\hat{\chi}(y_j)}{\hat{i}(y_j)} &=& 0. \nn\eea}

Let $F_{a, K}$ be a functional defined by (\ref{2.4}) for some $a, K\in\R$ large enough and let $\ep>0$ be
small enough such that $[-\ep, 0)$ contains no critical values of $F_{a, K}$. For $b$ large enough,
The normalized Morse series of $F_{a, K}$ in $ X^{-\ep}\setminus X^{-b}$
is defined, as usual, by
\be  M_a(t)=\sum_{q\ge 0,\;1\le j\le p} \dim C_{S^1,\;q}(F_{a, K}, \;S^1\cdot v_j)t^{q-d(K)},  \lb{2.11}\ee
where we denote by $\{S^1\cdot v_1, \ldots, S^1\cdot v_p\}$ the critical orbits of $F_{a, K}$ with critical
values less than $-\ep$. The Poincar\'e series of $H_{S^1, *}( X, X^{-\ep})$ is $t^{d(K)}Q_a(t)$, according
to Theorem 5.1 of \cite{LLW1}, if we set $Q_a(t)=\sum_{k\in \Z}{q_kt^k}$, then
$$   q_k=0 \qquad\qquad \forall\;k\in \mathring {I},  $$
where $I$ is an interval of $\Z$ such that $I \cap [i(\tau, y), i(\tau, y)+\nu(\tau, y)-1]=\emptyset$ for all
closed characteristics $(\tau,\, y)$ on $\Sigma$ with $\tau\ge aT$. Then by Section 6 of \cite{LLW1}, we have
$$  M_a(t)-\frac{1}{1-t^2}+Q_a(t) = (1+t)U_a(t),   $$
where $U_a(t)=\sum_{i\in \Z}{u_it^i}$ is a Laurent series with nonnegative coefficients.
If there is no closed characteristic with $\hat{i}=0$, then
\be   M(t)-\frac{1}{1-t^2}=(1+t)U(t),    \lb{2.12}\ee
where $M(t)$ denotes $M_a(t)$ as $a$ tends to infinity.

\setcounter{equation}{0}
\section{Multiplicity of closed characteristics on
star-shaped hypersurfaces}
In this section, we prove Theorems 1.1 based on Theorem 2.8, the index iteration theory developed by Y. Long
and his coworkers, and Theorem 1.2 of \cite{GHHM}.

The following theorem relates the Morse index defined in Section 2 to the Maslov-type index.

{\bf Theorem 3.1.} (cf. Theorem 2.1 of \cite{HuL1}) {\it Suppose $\Sg\in \H_{st}(2n)$ and
$(\tau,y)\in \T(\Sigma)$. Then we have
$$ i(y^m)\equiv i(m\tau,y)=i(y, m)-n,\quad \nu(y^m)\equiv\nu(m\tau, y)=\nu(y, m), \qquad \forall m\in\N, $$
where $i(y, m)$ and $\nu(y, m)$ are the Maslov-type index and nullity of $(m\tau,y)$ (cf. Section
5.4 of \cite{Lon2}). In particular, we have $\hat{i}(\tau,y)=\hat{i}(y,1)$, where $\hat{i}(\tau ,y)$ is
given in Section 2, $\hat{i}(y,1)$ is the mean Maslov-type index. Hence we denote it simply by
$\hat{i}(y)$.}

Now we prove Theorem 1.1. In the following, we fix $n=2$.

{\bf Proof of Theorem 1.1.} We prove Theorem 1.1 by contradiction. Now we assume
$\,^{\#}\T(\Sg)= 1$ for some $\Sg\in\H_{st}(4)$.

Denote the only one prime closed characteristic on $\Sg$ by $(\tau,\,y)$ with the corresponding associated
symplectic path $\ga\equiv \gamma_y:[0,\tau]\to\Sp(4)$. Then by Lemma 3.3 of \cite{HuL1} (cf. also Lemma
15.2.4 of \cite{Lon2}), there exist $P\in \Sp(4)$ and $M\in \Sp(2)$ such that
$\gamma(\tau)=P^{-1}(N_1(1,\,1)\diamond M)P$. By Theorem 2.8, we obtain $\hat{i}(y)>0$ and the following
identity
\be  \sum_{1\le m\le K(y)}\sum_{0\le l\le \nu(y^m)-1}
  \frac{(-1)^{i(y^{m})+l}k_l(y^{m})}{K(y)\hat{i}(y)}=\frac{1}{2}. \lb{3.1}\ee

If all the iterates of $(\tau, y)$ are non-degenerate, by Theorem 1.1 of \cite{HuL1}, there exist at least
two geometrically distinct closed characteristics. Thus to prove Theorem 1.1, it suffices to consider the
following four possible degenerate cases according to the classification of basic norm form decomposition
of $\gamma(\tau)$. In the following we use the notations from Definition 1.8.5 and Theorem 1.8.10 of
\cite{Lon2}:

{\bf Case 1.} {\it $\gamma(\tau)$ can be connected to
$\left(\begin{array}{cc}
        1 & 1 \\
        0 & 1 \\ \end{array}\right)\diamond \left(\begin{array}{cc}
                                                   -1 & b \\
                                                    0 & -1 \\ \end{array}\right)$
within $\Omega^0(\gamma(\tau))$, where $b=0$ or $\pm 1$.}

In this case, $i(y,1)$ is even by Theorems 8.1.4 and 8.1.5 of \cite{Lon2}. By Theorem 1.3 of
\cite{Lon1}, we have
\bea
i(y,m) &=& m(i(y, 1)+1)-1, \quad{\rm for}\;\; b=1;\nn\\
i(y,m)&=& m(i(y,1)+1)-1-\frac{1+(-1)^m}{2}, \quad{\rm for}\;\;b=0, -1. \nn\eea
By Theorem 3.1, we obtain
\bea
i(y^m) &=& m(i(y, 1)+1)-3, \quad{\rm for}\;\;b=1;  \lb{3.2}\\
i(y^m) &=& m(i(y,1)+1)-3-\frac{1+(-1)^m}{2}, \quad{\rm for}\;\;b=0, -1. \lb{3.3}\eea
By Proposition 2.6, we have $K(y)=2$. By (\ref{3.2})-(\ref{3.3}), we have
\be  \hat{i}(y)=i(y, 1)+1.  \lb{3.4}\ee
Thus $i(y, 1)\geq 0$ holds because $\hat{i}(y)>0$.

Now we proceed our proof in two subcases:

{\bf Subcase 1.1.} If $i(y, 1)= 0$, then $\hat{i}(y)=1$ by (\ref{3.4}). Noticing $K(y)=2$,
it follows from (\ref{3.1}) and (\ref{2.10}) that
$$  \frac{1+(-1)^{i(y^2)}(k_0(y^2)-k_1(y^2)+k_2(y^2))}{2}=\frac{1}{2},  $$
and then $k_0(y^2)-k_1(y^2)+k_2(y^2)=0$. Thus $k_0(y^2)=k_1(y^2)=k_2(y^2)=0$ by (v) of Remark 2.7. In this
subcase, by Proposition 2.4, we have $C_{S^1, q}(F_{a, K},S^1\cdot x^{2m})=0$ for all $q\in\Z$, and
$C_{S^1,q}(F_{a, K},S^1\cdot x^{2m-1})=\Q$ for $q=i(y^{2m-1})+d(K)$, $C_{S^1,q}(F_{a, K},S^1\cdot x^{2m-1})=0$ for
$q\neq i(y^{2m-1})+d(K)$, where $x$ is the critical point of $F_{a, K}$ corresponding to $y$,
$d(K)=4([KT/{2\pi}]+1)$, $m\in\N$. Note that $i(y^{2m-1})$ is even by (\ref{3.2})-(\ref{3.3}).

{\bf Subcase 1.2.} If $i(y, 1)\geq 2$, we have
\be   \frac{1+(-1)^{i(y^2)}(k_0(y^2)-k_1(y^2)+k_2(y^2))}{2(i(y,1)+1)}=\frac{1}{2}.   \lb{3.5}\ee
If $k_0(y^2)=1$ or $k_2(y^2)=1$, then by (i) and (v) of Remark 2.7, we have $k_0(y^2)-k_1(y^2)+k_2(y^2)=1$
and the left hand side of (\ref{3.5}) is less than $\frac{2}{2(i(y, 1)+1)}$ and we have
$\frac{2}{2(i(y, 1)+1)}\leq \frac{1}{3}$, which yields a contradiction with (\ref{3.5}). Hence
$k_0(y^2)=k_2(y^2)=0$ holds, and then $k_1(y^2)=i(y, 1)$ and $i(y^2)\in 2\Z+1$ must hold by (\ref{3.5}).
In this subcase, by Proposition 2.4, we have
\bea
C_{S^1,q}(F_{a, K},S^1\cdot x^{2m-1})=\Q, &&\quad{\rm for}\;\;q=i(y^{2m-1})+d(K),  \nn\\
C_{S^1,q}(F_{a, K},S^1\cdot x^{2m-1})=0, &&\quad{\rm for}\;\;q\neq i(y^{2m-1})+d(K), \nn\\
C_{S^1,q}(F_{a, K},S^1\cdot x^{2m})=\Q, &&\quad{\rm for}\;\;q=i(y^{2m})+1+d(K),  \nn\\
C_{S^1, q}(F_{a,K},S^1\cdot x^{2m})=0, &&\quad{\rm for}\;\;q\neq i(y^{2m})+1+d(K), \nn\eea
where $x$ is the critical point of $F_{a, K}$ corresponding to $y$, $d(K)=4([KT/{2\pi}]+1)$, and $m\in\N$.
Note that $i(y^{2m})$ is odd in this subcase.

In any subcase, we have $C_{S^1, q}(F_{a, K},S^1\cdot x^{m})=0$ when $q$ is odd, $\forall~m\in \N$, then we have
that the normalized Morse series $M(t)$ of $F_{a, K}$ (cf. (\ref{2.11}) and (\ref{2.12})) have only even
power terms. We apply (\ref{2.12}) to obtain that $U(t)\equiv 0$, and then
\be   M(t)-\sum_{m=0}^\infty{t^{2m}}=0.  \lb{3.6}\ee
For Subcase 1.1, by (\ref{2.11}) and (\ref{3.2})-(\ref{3.3}), we have
\be   M(t)=\sum_{m=1}^\infty t^{i(y^{2m-1})}=\sum_{m=1}^\infty t^{2m-4}.  \lb{3.7}\ee
Comparing (\ref{3.7}) and (\ref{3.8}), we get a contradiction!

For Subcase 1.2, by (\ref{2.11}), (\ref{3.2})-(\ref{3.3}) and $k_1(y^2)=i(y, 1)$ which we obtain in
Subcase 1.2, we have
\bea  M(t)
&=& t^{i(y)}+k_1(y^2)t^{i(y^2)+1}+t^{i(y^3)}+\cdots   \nn\\
&=& t^{i(y)}+i(y,1)t^{i(y^2)+1}+t^{i(y^3)}+\cdots.    \lb{3.8}\eea
From (\ref{3.6}) and (\ref{3.8}), and the assumption that $i(y, 1)\geq 2$ in subcase 1.2, we get a
contradiction!

{\bf Case 2.} {\it $\gamma(\tau)$ can be connected to
$\left(\begin{array}{cc}
        1 & 1 \\
        0 & 1 \\ \end{array}\right)\diamond \left(\begin{array}{cc}
                                                  \cos{\theta} & -\sin{\theta} \\
                                                  \sin{\theta} & \cos{\theta} \\ \end{array}\right)$
within $\Omega^0(\gamma(\tau))$, $\theta\in(0, \pi)\cup (\pi, 2\pi)$ and $\theta/\pi\in {\bf Q}$.}

In this case, $i(y, 1)$ is even by Theorems 8.1.4 and 8.1.7 of \cite{Lon2}. By Theorem 1.3 of \cite{Lon1},
we have
$$ i(y,m)=mi(y,1)+2E(\frac{m\theta}{2\pi})-2,  $$
where $E(a)=\min{\{k\in\Z\mid k\geq a\}}$ for any $a\in\R$. By Theorem 3.1, we obtain
\be  i(y^m)=mi(y,1)+2E(\frac{m\theta}{2\pi})-4.  \lb{3.9}\ee
Then
\be  \hat{i}(y)=i(y,1)+\frac{\theta}{\pi}.  \lb{3.10}\ee
Since $\hat{i}(y)>0$ and the fact that $i(y, 1)$ is even, we get $i(y, 1)\geq 0$.

If $i(y, 1)\geq 2$. Then $\hat{i}(y)>2$ by (\ref{3.10}). Note that $k_0(y^m)=1$ and $k_l(y^m)=0$ for
$1\leq m\leq K(y)-1$ and $l\neq 0$ by (\ref{2.10}). Now (\ref{3.1})
becomes the following
\be \frac{K(y)-1+k_0(y^{K(y)})-k_1(y^{K(y)})+k_2(y^{K(y)})}{K(y)\hat{i}(y)}=\frac{1}{2}.  \lb{3.11}\ee
Then
$$  \frac{K(y)-1+k_0(y^{K(y)})-k_1(y^{K(y)})+k_2(y^{K(y)})}{K(y)}=\frac{\hat{i}(y)}{2}>1,  $$
which contradicts to (i) and (v) of Remark 2.7.

Hence, in the rest of Case 2 we suppose $i(y, 1)= 0$, then by (\ref{3.10}) we have
\be  \hat{i}(y)=\theta/\pi<2.  \lb{3.12}\ee

We claim that $k_0(y^{K(y)})=k_2(y^{K(y)})=0$. Otherwise, by Remark 2.7 (iv), one of $k_0(y^{K(y)})$
and $k_2(y^{K(y)})$ is equal to 1 and $k_1(y^{K(y)})=0$, then
$$  K(y)-1+k_0(y^{K(y)})-k_1(y^{K(y)})+k_2(y^{K(y)})=K(y)  $$
holds. Together with (\ref{3.11}) it yields $\hat{i}(y)=2$, which contradicts to (\ref{3.12}). Hence,
(\ref{3.11}) becomes
\be \frac{2(K(y)-1-k_1(y^{K(y)}))}{K(y)} = \frac{\theta}{\pi}.   \lb{3.13}\ee
The normalized Morse series $M(t)$ of $F_{a, K}$ has the following form:
\be  M(t) = t^{i(y)}+t^{i(y^2)}+\cdots+t^{i(y^{K(y)-1})}+k_1(y^{K(y)})t^{i(y^{K(y)})+1}+\cdots.  \lb{3.14}\ee
We apply (\ref{2.12}) to obtain
\be  M(t)-\sum_{m=0}^\infty{t^{2m}} = (1+t)U(t),   \lb{3.15}\ee
where $U(t)$ is a Laurent series with nonnegative coefficients.

Noticing that $i(y)=i(y,1)-2=-2$, $i(y^m)\geq i(y)=-2$ and it is even for all $m\in{\bf N}$ by (\ref{3.9}).
Then from (\ref{3.14}) and (\ref{3.15}) we can suppose that $U(t)$ has the form:
$$  U(t)=a_{-2}t^{-2}+\sum_{i\geq-1}{a_i t^i},   $$
and $a_{-2}\in {\bf N}$. Thus the coefficient of $t^{-1}$ in the right hand side of (\ref{3.15}) is
nonzero, and then we get $i(y^{K(y)})=-2$, and then $i(y)=\cdots=i(y^{K(y)-1})=-2$. Comparing the
coefficient of $t^{-2}$ in the two sides of (\ref{3.15}), we obtain $a_{-2}=K(y)-1$. Then comparing the
coefficient of $t^{-1}$ in the two sides of (\ref{3.15}), we obtain $k_1(y^{K(y)})= a_{-2}+a_{-1}\geq K(y)-1$,
which contradicts to (\ref{3.13}).

{\bf Case 3.} {\it $\gamma(\tau)$ can be connected to
$\left(\begin{array}{cc}
        1 & 1 \\
        0 & 1 \\ \end{array}\right)\diamond\left(\begin{array}{cc}
                                                  1 & b \\
                                                  0 & 1 \\ \end{array}\right)$
within $\Omega^0(\gamma(\tau))$, where $b=0$ or $1$.}

In this case, $i(y, 1)$ is even by Theorem 8.1.4 of \cite{Lon2}. By Theorem 1.3 of \cite{Lon1}, we have
$$ i(y,m) = m(i(y, 1)+2)-2. $$
By Theorem 3.1, we obtain
\be  i(y^m) = m(i(y,1)+2)-4.   \lb{3.16}\ee
By Proposition 2.6, we have $K(y)=1$. By (\ref{3.16}), we have
\be  \hat{i}(y)=i(y, 1)+2.   \lb{3.17}\ee
Since $\hat{i}(y)>0$ and the fact $i(y,1)$ is even, then we have $i(y,1)\geq 0$.

In this case, (\ref{3.1}) becomes:
\be  \frac{(-1)^{i(y)}(k_0(y)-k_1(y)+k_2(y))}{i(y, 1)+2}=\frac{1}{2}.  \lb{3.18}\ee
Since $i(y)$ is even, by (v) of Remark 2.7, we obtain one of $k_0(y)$ and $k_2(y)$ is equal to 1 and $k_1(y)=0$.
Then $k_0(y)-k_1(y)+k_2(y)=1$ and $i(y, 1)=0$ by (\ref{3.18}), thus $\hat{i}(y)=2$ by (\ref{3.17}).

In the following, we prove that the prime closed characteristic $(\tau, y)$ satisfies (LF3)-(ii) on Page 10 of
\cite{GHHM}, then by Theorem 1.2 of \cite{GHHM}, we obtain the existence of infinitely many geometrically distinct
prime closed characteristics which contradicts to our assumption. This completes the proof in Case 3.

Let $\eta$ be the contact form on $\Sigma$ arising from $\Sigma$ being a $C^3$ compact hypersurface in $\R^{4}$
strictly star-shaped with respect to the origin. The key property of $\eta$ we need is that the contact homology
$HC_*(\eta)$ is concentrated in even dimensions $2, 4, \cdots$ and that $HC_*(\eta)=\Q$ in these dimensions
(cf. \cite{Bou1} and references therein). Here the indices are evaluated with respect to a global symplectic
trivialization of $ker\eta$ on $\Sigma$; such a trivialization is unique up to homotopy. Thus by (LF1)
on Page 9 of \cite{GHHM}, we get the local contact homology $HC_2(y)\neq 0$ (cf. \cite{HrM1} for the definition
of the local contact homology), noticing that $\hat{i}(y)=2$ and $(\tau, y)$ is degenerate, then $(\tau, y)$
satisfies (LF3)-(ii) on Page 10 of \cite{GHHM}.

{\bf Case 4.} {\it $\gamma(\tau)$ can be connected to
$\left(\begin{array}{cc}
        1 & 1 \\
        0 & 1 \\ \end{array}\right)\diamond\left(\begin{array}{cc}
                                                  1 & -1 \\
                                                  0 & 1 \\ \end{array}\right)$
within $\Omega^0(\gamma(\tau))$.}

In this case, $i(y, 1)$ is odd by Theorem 8.1.4 of \cite{Lon2}. By Theorem 1.3 of \cite{Lon1}, we have
$$  i(y,m) = m(i(y, 1)+1)-1.  $$
By Theorem 3.1, we obtain
\be  i(y^m)=m(i(y,1)+1)-3.   \lb{3.19}\ee
By Proposition 2.6, we have $K(y)=1$. By (\ref{3.19}), we have
\be  \hat{i}(y)=i(y, 1)+1.   \lb{3.20}\ee
Since $\hat{i}(y)>0$ and the fact that $i(y, 1)$ is odd, then $i(y, 1)\geq 1$. From (\ref{3.1}), we have
\be  \frac{(-1)^{i(y)}(k_0(y)-k_1(y)+k_2(y))}{i(y, 1)+1}=\frac{1}{2}.  \lb{3.21}\ee
Since $i(y)$ is odd, then $k_0(y)=k_2(y)=0$ by Remark 2.7 (iv) and (\ref{3.21}) becomes:
\be  \frac{k_1(y)}{i(y,1)+1}=\frac{1}{2}.  \lb{3.22}\ee
Note that $i(y^m)$ is always odd for all $m\in \N$ by (\ref{3.19}), then the normalized Morse series
of $F_{a, K}$ is $M(t)=\sum_{m\geq1}k_1(y)t^{i(y^m)+1}$ which has only even power terms. We apply (\ref{2.12})
to obtain $U(t)\equiv 0$ and
$$ \sum_{m\geq1}k_1(y)t^{i(y^m)+1}-\sum_{m=0}^\infty{t^{2m}}=0. $$
Then $k_1(y)=1$ must hold. Together with (\ref{3.20}) and (\ref{3.22}), it yields $\hat{i}(y)=i(y, 1)+1=2$.
Then by the same argument in the last step of the proof for Case 3, we complete the proof for Case 4.

The proof of Theorem 1.1 is completed.

{\bf Acknowledgements.} The authors would like to sincerely thank the anonymous referee for his
careful reading of the manuscript and valuable comments.

\bibliographystyle{abbrv}

\end{document}